\documentclass[12pt,a4paper,reqno]{amsart}
\usepackage[english]{babel}
\usepackage{amssymb,latexsym,amsfonts,amsthm,upref,amsmath}
\usepackage[foot]{amsaddr}
\usepackage[margin=1in]{geometry} 
\usepackage[dvipsnames,x11names]{xcolor}
 \definecolor{myblue}{HTML}{003399}
\usepackage{hyperref}
\hypersetup{colorlinks,citecolor=myblue,filecolor=black,linkcolor=myblue,urlcolor=myblue}
\usepackage{enumerate}  
\usepackage{tikz}
\usepackage{float}
\usepackage{cleveref}
\usepackage{mathtools}
\usepackage{cite}
\usepackage{filecontents}
\usepackage{enumitem}
\makeatletter
\newcommand{\leqnomode}{\tagsleft@true}
\newcommand{\reqnomode}{\tagsleft@false}

\makeatother
\newtheorem*{thm*}{Theorem}
\newtheorem*{lem*}{Lemma}
\newtheoremstyle{prim}{}{}{\normalfont}{}{\bfseries}{.}{ }{}
\newtheoremstyle{stil}{}{}{\slshape}{}{\bfseries}{.}{ }{}
\theoremstyle{stil}
\newtheorem{thm}{Theorem}[section]
\newtheoremstyle{defi}{}{}{}{}{\bfseries}{.}{ }{}
\theoremstyle{defi}

\theoremstyle{defi}
\newtheorem{rem}[thm]{Remark}
\theoremstyle{stil}
\newtheorem*{mthm*}{Main Theorem}
\newtheorem*{kor*}{Corollary}
\newtheorem{pro}[thm]{Proposition}
\theoremstyle{stil}
\newtheorem{lem}[thm]{Lemma}
\theoremstyle{stil}
\newtheorem{kor}[thm]{Corollary}
\theoremstyle{prim}

\newenvironment{prf}{\noindent \textit{Proof.}}{\null\hfill$\qed$\hskip
2mm\vskip 2mm}

\newcommand{\dyo}{{\rm DY}_0(\mathfrak{gl}_{m|n})}
\newcommand{\dyoc}{ \wtld{{\rm DY}}_0(\mathfrak{gl}_{m|n})}
\newcommand{\dyc}{ \wtld{{\rm DY}}(\mathfrak{gl}_{m|n})}

\newcommand{\dyg}{ {\rm DY}(\mathfrak{gl}_{m|n})}

\newcommand{\Y}{ {\rm Y}(\mathfrak{gl}_{m|n})}
\newcommand{\Yd}{ {\rm Y}^+(\mathfrak{gl}_{m|n})}
\newcommand{\Ydt}{ \widetilde{\rm Y}^+(\mathfrak{gl}_{m|n})}

\newcommand{\U}{ {\rm U}}

\newcommand{\B}{ {\rm B}(\mathfrak{gl}_{m|n})}
\newcommand{\Bd}{ {\rm B}^+(\mathfrak{gl}_{m|n})}
\newcommand{\DB}{ {\rm DB} (\mathfrak{gl}_{m|n})}

\newcommand{\DBB}{ \mathcal{  DB}  (\mathfrak{gl}_{m|n})}
\newcommand{\BB}{ \mathcal{  B}  (\mathfrak{gl}_{m|n})}
\newcommand{\BBd}{ \mathcal{  B}^+  (\mathfrak{gl}_{m|n})}

\newcommand{\Vccgl}{\wtld{\rm V}_{\hspace{-1pt}c}(\gl_{m|n})}

\newcommand{\R}{ {\overline{R}}}

\newcommand{\vac}{\mathop{\mathrm{\boldsymbol{1}}}}

\newcommand{\gr}{\mathop{\mathrm{gr}}}

\newcommand{\glmn}{\mathfrak{gl}_{m|n}}
\newcommand{\glmnht}{\widehat{\mathfrak{gl}}_{m|n}}
\newcommand{\gl}{\mathfrak{gl}}

\newcommand{\CC}{\mathbb{C}}

\newcommand{\ZZ}{\mathbb{Z}}

\newcommand{\wtld}{\widetilde}

\newcommand{\wvr}{\overline}

\newcommand{\ot}{\otimes}
\newcommand{\ts}{\hspace{1pt}}
\newcommand{\qdet}{ \mathop{\rm qdet} }

\newcommand{\ndo}{\mathop{\mathrm{End}}}

\newcommand{\diag}{\mathop{\mathrm{diag}}}

\newcommand{\fand}{\quad\text{and}\quad}
\newcommand{\Fand}{\qquad\text{and}\qquad}

\newcommand{\non}{\nonumber}
\newcommand{\beq}{\begin{equation}}
\newcommand{\eeq}{\end{equation}}
\newcommand{\ben}{\begin{equation*}}
\newcommand{\een}{\end{equation*}}

\makeatletter
\def\smalloverbrace#1{\mathop{\vbox{\m@th\ialign{##\crcr\noalign{\kern3\p@}%
  \tiny\downbracefill\crcr\noalign{\kern3\p@\nointerlineskip}%
  $\hfil\displaystyle{#1}\hfil$\crcr}}}\limits}
\makeatother

\makeatletter
\def\smallunderbrace#1{\mathop{\vtop{\m@th\ialign{##\crcr
   $\hfil\displaystyle{#1}\hfil$\crcr
   \noalign{\kern3\p@\nointerlineskip}%
   \tiny\upbracefill\crcr\noalign{\kern3\p@}}}}\limits}
\makeatother

\begin{document}

\title{Double Yangian and  reflection algebras of the Lie superalgebra $\mathfrak{gl}_{m|n}$}

\author{Lucia Bagnoli}
\author{Slaven Ko\v{z}i\'{c}}
\address[L. Bagnoli and S. Ko\v{z}i\'{c}]{Department of Mathematics, Faculty of Science, University of Zagreb,  Bijeni\v{c}ka cesta 30, 10\,000 Zagreb, Croatia}
\email{lucia.bagnoli@math.hr}
\email{kslaven@math.hr}

 \begin{abstract}
We study the  double Yangian  associated with the Lie superalgebra $\mathfrak{gl}_{m|n}$.
Our main focus is on establishing the Poincar\'{e}--Birkhoff--Witt Theorem for the double Yangian and   constructing   its central elements in the form of coefficients of the quantum contraction. Next, as an application, we introduce reflection algebras,     certain left coideal subalgebras of the level 0 double Yangian,   and find their presentations by generators and relations.
\end{abstract}

 \maketitle

\allowdisplaybreaks

\section{Introduction}\label{intro}
\numberwithin{equation}{section}

The  {\em Yangian} $\Y$ for the  general linear Lie superalgebra $\mathfrak{gl}_{m|n}$ was introduced by Nazarov \cite{Naz} via an $R$-matrix presentation. It can be viewed as a deformation of the universal enveloping algebra $\U(\glmn[t])$. Since its introduction, its structure was extensively studied, in particular, due to its close connections with various areas of mathematical physics, such as Calogero-Sutherland systems \cite{AK,JWW}, non-linear super-Schr\"{o}dinger equation \cite{CR} and superstrings on $\text{AdS}_5\times \text{S}^5$ \cite{HY}. The super Yangian possesses two distinct families of central elements, established in \cite{Naz}, which consist of coefficients of certain formal power series, {\em quantum Berezinian} $b(u)$ and {\em quantum contraction} $z(u)$. These series can be regarded as  super analogues of the   quantum determinant  $\qdet T(u)$ for the ordinary Yangian $ {\rm Y}(\mathfrak{gl}_{N})$ and the series $\qdet T(u-1)/\qdet T(u)$, respectively; see, e.g.,  \cite[Ch. 1]{Mol}. 
As with their even counterparts, the  coefficients of $b(u)$ and $z(u)$ generate the  entire centre of  $\Y$, which was conjectured by Nazarov \cite{Naz} and proved by Gow \cite{Gow}. 

In this paper,   we consider the   {\em double Yangian} $\dyg$ for the  general linear Lie superalgebra $\mathfrak{gl}_{m|n}$ given by an $R$-matrix presentation, which is based on the definition of Zhang \cite{Z}.
Our main goal is to establish the Poincar\'{e}--Birkhoff--Witt Theorem for $\dyg$.
In the even case, 
the Poincar\'{e}--Birkhoff--Witt Theorem for  ${\rm DY}(\mathfrak{gl}_N)$ goes back to the papers by   Jing, Molev, Yang and the second author \cite{JKMY}, Nazarov  \cite{Naz3} and also the paper by   Wendlandt  \cite{W}, where it was proved for the double Yangian of an arbitrary finite-dimensional or simply laced affine Kac--Moody  Lie algebra. The proof of this theorem in the super case, which we give in Section  \ref{section1}, relies, in particular, on the ideas of Etingof and Kazhdan \cite{EK3, EK} and Nazarov \cite{Naz2,Naz3}. 
Furthermore, it employs certain preliminary results on the {\em dual Yangian} $\Yd$ and the {\em level  0  double Yangian} $\dyo$, so we study these algebras before proceeding towards $\dyg$.

Our second goal in Section  \ref{section1} is to investigate possible generalizations of the quantum contraction $z(u)$ for the Yangian  $\Y$ to the  double Yangian $\dyg$. 
This research direction is partially motivated by the fact that the series $z(u)$ has  not been   widely studied in the literature, in contrast with the quantum Berezinian $b(u)$; see, e.g., the papers \cite{Gow0,LM,MR2,T}. We introduce the dual Yangian analogue $z^+(u)$ of the quantum contraction and we show that, as with $z(u)$, its coefficients are algebraically independent elements of the centre of the suitably completed algebra  $\dyg$.

In Section \ref{section2}, we consider the so-called {\em reflection algebras}. The algebras associated with the  {\em reflection equation} were originally  introduced by Sklyanin \cite{S} to describe integrable systems with the boundary conditions; see also \cite{GS,KS,KJC,MRS}   for more information on such algebras and their applications.
A  distinct class of reflection algebras, which are left coideal subalgebras in the Yangian $ {\rm Y}(\mathfrak{gl}_{N})$, was studied by Molev and Ragoucy \cite{MR}.
Moreover, its connection with Etingof--Kazhdan's quantum affine vertex algebras \cite{EK} was investigated by the second author 
\cite{K}. In this section, generalizing the approach from \cite{MR,K} to the super case, we introduce {\em double reflection algebras} $\DB$, which are left coideal subalgebras of  (a suitable completion of)  $\dyo$. As an application of the  Poincar\'{e}--Birkhoff--Witt Theorem for $\dyg$, which implies that there exists an isomorphism of $\ZZ_2$-graded algebras
$$
  \textstyle\gr_2 \dyg \cong \U(\glmnht),
$$
we obtain the isomorphism
$$
\textstyle\gr_2 \DB\cong \U(\mathfrak{gl}_{m|n}[t,t^{-1}]^{\sigma} ).
$$
Here  $\gr_2 A$ stands for the corresponding graded algebra of the algebra $A$ with respect to a certain degree operator and $\mathfrak{gl}_{m|n}[t,t^{-1}]^{\sigma}$ is a   subalgebra of the loop Lie superalgebra $\mathfrak{gl}_{m|n}[t,t^{-1}] $ which depends on the choice of   involutive automorphism $\sigma$ of $\mathfrak{gl}_{m|n}$.
The main result of   Section \ref{section2} is an explicit presentation of the subalgebra $\DB\subset \dyo$. More specifically, we show that it can be defined as an algebra in given generators subject to the family of three reflection relations and two unitarity constraints.

\section{Double Yangian for the Lie superalgebra \texorpdfstring{$\glmn$}{glm|n}}\label{section1}

In this section, we recall some  properties of the Yangian $\Y$. Next, we study the dual Yangian    $\Yd$, the   level  0  double Yangian  $\dyo$ and the (centrally extended) double Yangian $\dyg$ with $m\neq n$ and establish our main result, the Poincar\'{e}--Birkhoff--Witt Theorem for $\dyg$. Our definition of the double Yangian closely follows Zhang \cite{Z}, but we use a different normalization of the Yang $R$-matrix which governs its defining relations.  Finally, we introduce the dual Yangian analogue of the quantum contraction and study its properties.

\subsection{Yangian for \texorpdfstring{$\glmn$}{glm|n}}\label{subsec0101}

Consider the Lie superalgebra 
$\glmnht=\glmn\ot\CC[t,t^{-1}]\oplus\CC K .$
Its supercommutation relations are given by
\begin{align}
[e_{ij}(r),e_{kl}(s)]
=&\, 
\delta_{kj}\ts  e_{il}(r+s)
-  \delta_{il}\ts e_{kj}(r+s) (-1)^{(\bar{i}+\bar{j})(\bar{k}+\bar{l})}\non\\
&+K\left(  \delta_{kj}\ts \delta_{il}(-1)^{\bar{i}} -\frac{\delta_{ij}\ts \delta_{kl}}{m-n}\ts (-1)^{\bar{i}+\bar{k}}\right)\ts r\ts \delta_{r+s\ts 0}\label{glmn}
\end{align}
for $m\neq n$,
where 
the element $K$ is even and central,  
$e_{ij}\in \glmn$ are  matrix units and $e_{ij}(r)=e_{ij}\ot t^r$. 
The parity of
the element
$e_{ij}(r)$ is $\bar{i}+\bar{j}$, where  $\bar{i}=0$ for $i=1,\ldots ,m$ and $\bar{i}=1$ for $i=m+1,\ldots ,m+n$.
Note that we can rescale the central element by $K=(n-m)K'$ so that the above relations apply to the case $m=n$ as well.

We  follow Nazarov \cite{Naz} to define the   Yangian for the general linear Lie superalgebra. The {\em Yangian} $\Y$ is the $\ZZ_2$-graded unital associative algebra with generators $t_{ij}^{(r)}$, where $  i,j =1,\ldots , m+n$ and $r=1,2,\ldots ,$  subject to the defining relations
\beq\label{Y}
[t_{ij}^{(r)},t_{kl}^{(s)}]
=(-1)^{\bar{i}\bar{j}+\bar{i}\bar{k}+\bar{j}\bar{k}}
\sum_{a=1}^{\min\left\{r,s\right\}}\left(
t_{kj}^{(a-1)}\ts t_{il}^{(r+s-a)}-t_{kj}^{(r+s-a)}\ts  t_{il}^{(a-1)}\right),
\eeq
where the square brackets denote the supercommutator and $t_{ij}^{(0)}=\delta_{ij}$. The parity of
the element
$t_{ij}^{(r)}$ is $\bar{i}+\bar{j}$. Relations \eqref{Y} can be  expressed in terms of   formal power series
$$
t_{ij}(u)=\delta_{ij}+\sum_{r\geqslant 1} t_{ij}^{(r)} u^{-r}
$$
as
\beq\label{Y2}
[t_{ij} (u),t_{kl} (v)]
=\frac{(-1)^{\bar{i}\bar{j}+\bar{i}\bar{k}+\bar{j}\bar{k}}}{u-v}
\left(t_{kj}(u)\ts t_{il}(v)-t_{kj}(v)\ts t_{il}(u)\right).
\eeq

Consider the rational $R$-matrix $R(u)\in\ndo\CC^{m|n}\ot \ndo\CC^{m|n} [u^{-1}]$ given by
\beq\label{Rmatrix}
R(u)=1-Pu^{-1},\quad\text{where}\quad P=\sum_{i,j=1}^{m+n} e_{ij}\ot e_{ji}\ts (-1)^{\bar{j}}
\eeq
and $1$ is the identity. The $R$-matrix satisfies the {\em Yang--Baxter equation}
\beq\label{ybe}
R_{12}(u_1-u_2)\ts R_{13}(u_1-u_3)\ts R_{23}(u_2 -u_3)=
R_{23}(u_2 -u_3)\ts R_{13}(u_1-u_3)\ts R_{12}(u_1-u_2).
\eeq
By using the $R$-matrix $R(u)$, one can write the defining relations for the super Yangian    in the so-called $RTT$-form as follows. 
Organize the series $t_{ij}(u)$ into the matrix
\beq\label{teu}
T(u)=\sum_{i,j=1}^{m+n} (-1)^{\bar{i}\bar{j}+\bar{j}} e_{ij}\ot t_{ij}(u)
 .
\eeq
Relations \eqref{Y}   are then expressed as the   identity of formal power series in the variables $u$ and $v$ such that their coefficients belong to $\ndo\CC^{m|n}\ot\ndo\CC^{m|n}\ot\Y $,
\beq\label{YRTT}
R(u-v)\ts T_1(u)\ts T_2(v)
=T_2(v)\ts T_1(u)\ts R(u-v).
\eeq
Note that in \eqref{YRTT}  we used the subscripts to  indicate the tensor factors on which the corresponding matrices are applied, so that we have
$$
T_1(u)=\sum_{i,j=1}^{m+n} (-1)^{\bar{i}\bar{j}+\bar{j}} e_{ij}\ot 1\ot t_{ij}(u)\Fand
T_2(v)=\sum_{i,j=1}^{m+n} (-1)^{\bar{i}\bar{j}+\bar{j}} 1\ot e_{ij}\ot t_{ij}(v).
$$

We now recall the Poincar\'{e}--Birkhoff--Witt Theorem for $\Y$  which was proved by Gow \cite{Gow}. First, introduce two different ascending filtrations on the   Yangian   by
\beq\label{deg21}
\deg_1 t_{ij}^{(r)} = r\fand \deg_2 t_{ij}^{(r)}=r-1.
\eeq
Denote by $\gr_1 \Y$ and $\gr_2 \Y$ the corresponding graded algebras.
We use the notation $\hat{t}_{ij}^{(r)}$ (resp. $\bar{t}_{ij}^{(r)}$) for the images of  generators in the respective components of the graded algebra    $\gr_1 \Y$ (resp. $\gr_2 \Y$). By \eqref{Y}, $\gr_1 \Y$ is supercommutative, i.e. we have
$
[\hat{t}_{ij}^{(r)},\hat{t}_{kl}^{(s)}]=0.
$
On the other hand,  taking the images of the generators in the respective components of   $\gr_2 \Y$ in \eqref{Y}, one   checks   that the assignments
\beq\label{assignments1}
e_{ij}(r-1)\mapsto (-1)^{\bar{i}}\ts \bar{t}_{ij}^{(r)},
\eeq
where  $i,j=1,\ldots ,m+n$ and $r\geqslant 1$,
define a homomorphism 
\beq\label{mapY}
\U(\glmn[t])\to \textstyle\gr_2\Y.
\eeq
Let us recall \cite[Thm. 1]{Gow}, which implies that the map   \eqref{mapY} is an algebra isomorphism.

\begin{thm}\label{pbw}
{\rm (1)}\label{xcvbnv}  Fix some ordering on the generators $t_{ij}^{(r)}$, where $i,j=1,\ldots ,m+n$ and $r\geqslant 1$. 
Then the ordered monomials in  $t_{ij}^{(r)}$, with at most power $1$ for odd generators, form a basis of the Yangian $\Y$.  
 \vspace{1pt}

\noindent{\rm (2)} The graded algebra $\gr_1 \Y$ is supercommutative and the corresponding images of the ordered monomials from   assertion  {\rm (\hyperref[xcvbnv]{1})} form its basis.
\end{thm}

Consider the transposition $ 
\tau\colon e_{ij}\mapsto (-1)^{\bar{i}\bar{j}+\bar{i}}e_{ji}
 $  on
$\ndo \CC^{m|n}$.
In  \cite{Naz}, Nazarov introduced a certain family of central elements of the Yangian. They were defined as coefficients of the {\em quantum contraction}, a  power series $z (u)\in\Y[[u^{-1}]] $ which is uniquely determined by the identity
\beq\label{qcontra}
P^{\tau_2}\ts  T_1(u+m-n) \left(T_2(u)^{-1}\right)^\tau = P^{\tau_2}\ot  z (u).
\eeq
Here $\tau_2$ indicates that the transposition $\tau$ is applied on the second tensor factor. 
Finally, we recall that the quantum analogue of the Liouville theorem    \cite[Thm. 2]{Naz}, along with the fact that the coefficients of the quantum Berezinian  generate the centre of the Yangian   \cite[Thm. 4]{Gow}, implies  that  the coefficients of $z (u)$ also generate the 
centre of  $\Y$.

\subsection{Dual Yangian for \texorpdfstring{$\glmn$}{glm|n}}\label{subsec0102}

Define the {\em dual Yangian} $\Yd$ as the $\ZZ_2$-graded unital associative algebra with generators $t_{ij}^{(-r)}$, where $  i,j=1,\ldots , m+n$ and $r=1,2,\ldots ,$ subject to the  defining relations 
\begin{align}
[t_{ij}^{(-r)},t_{kl}^{(-s)}]
=&\,(-1)^{\bar{i}\bar{j}+\bar{i}\bar{k}+\bar{j}\bar{k}}
\Big(
\delta_{kj}\ts t_{il}^{(-r-s)}-\delta_{il}\ts t_{kj}^{(-r-s)}
\Big.\non\\
\Big.
&+\sum_{a=1}^{\min\left\{r,s\right\}}
\left( t_{kj}^{(-r-s+a-1)} t_{il}^{(-a)}-t_{kj}^{(-a)} t_{il}^{(-r-s+a-1)}\right)\Big).\label{Yd}
\end{align}
As before, the square brackets   denote the supercommutator and the parity of   $t_{ij}^{(-r)}$
is $\bar{i}+\bar{j}$. Defining relations  \eqref{Yd} can be   expressed in terms of formal power series
$$
t_{ij}^+ (u)=\delta_{ij}-\sum_{r\geqslant 1} t_{ij}^{(-r)} u^{r-1}
$$
as
$$
[t_{ij}^+ (u),t_{kl}^+ (v)]
=\frac{(-1)^{\bar{i}\bar{j}+\bar{i}\bar{k}+\bar{j}\bar{k}}}{u-v}
\left(t_{kj}^+(u) \ts t_{il}^+(v)-t_{kj}^+(v)\ts t_{il}^+(u)\right).
$$
Moreover, as with the Yangian, one can use the $R$-matrix \eqref{Rmatrix} to  write defining relations \eqref{Yd} in the $RTT$-form, i.e. as the     identity of formal power series in the variables $u$ and $v$ such that their coefficients belong to $\ndo\CC^{m|n}\ot\ndo\CC^{m|n}\ot\Yd  $,
\beq\label{YdRTT}
R(u-v)\ts T_1^+ (u)\ts T_2^+ (v)
=T_2^+ (v)\ts T_1^+ (u)\ts R(u-v),
\eeq
 where  $T^+(u)$ is given by
\beq\label{teud}
T^+ (u)=\sum_{i,j=1}^{m+n} (-1)^{\bar{i}\bar{j}+\bar{j}} e_{ij}\ot t^+_{ij}(u)
 .
\eeq

To prove  the Poincar\'{e}--Birkhoff--Witt Theorem for the dual Yangian we   follow 
the approach of Nazarov   \cite[Sect. 11,  12]{Naz3}, where such a result was proved in the even case.
Introduce the ascending filtration on $\Yd$ by
\beq\label{deg22}
 \deg_2 t_{ij}^{(-r)}=-r.
\eeq
Denote by $\gr_2 \Yd$ the corresponding graded algebra.
We use the notation   $\bar{t}_{ij}^{(-r)}$  for the images of the generators    in the respective components of    $\gr_2 \Yd$. Defining relations \eqref{Yd} imply the identities in $\gr_2 \Yd$,
$$
[\bar{t}_{ij}^{(-r)},\bar{t}_{kl}^{(-s)}]
= 
(-1)^{ \bar{i} \bar{j} + \bar{i}\bar{k} + \bar{j}\bar{k}}
\left(
\delta_{kj}\ts \bar{t}_{il}^{(-r-s)}-\delta_{il} \ts\bar{t}_{kj}^{(-r-s)}\right) .
$$
 Therefore, by commutation  relations  \eqref{glmn} for the Lie superalgebra $\glmnht$,   the assignments
\beq\label{assignments2}
e_{ij}(-r)\mapsto (-1)^{\bar{i}}\ts \bar{t}_{ij}^{(-r)},
\eeq
where
$ i,j=1,\ldots ,m+n$ and $r\geqslant 1$,
define a homomorphism 
\beq\label{mapYdx}
\U(t^{-1}\glmn[t^{-1}])\to \textstyle\gr_2\Yd.
\eeq

 The aforementioned proof in \cite{Naz3} employs a certain bilinear pairing  motivated by  \cite[Sect. 2]{RTF}. In the super case, such a pairing $\left<\cdot,\cdot\right>\colon\Y\times\Yd\to \CC$     is defined   so that  the corresponding linear map   $\Y\otimes\Yd\to \CC$ satisfies 
$$
T_1(u_1)\ldots T_k(u_k)\ts T_{k+1}^+(v_1)\ldots T_{k+l}^+(v_l)
\mapsto
\prod_{i=1,\ldots, k }^{\longrightarrow}
\prod_{j=1,\ldots, l }^{\longleftarrow}
R_{i\ts j+k}(u_i -v_j)
$$
in 
$
(\ndo\CC^{m|n})^{\ot k} \ot (\ndo\CC^{m|n})^{\ot l}  \ot\Y\ot \Yd
$
for all integers $k,l\geqslant 0$, where the arrows indicate the order of  factors. In particular, we have $\left<1,1\right>=1$. The fact that the pairing is well-defined is easily proved by using the Yang--Baxter equation \eqref{ybe} and   the defining relations  in the $RTT$-form, \eqref{YRTT} and \eqref{YdRTT}. The following property of the pairing
can be verified by arguing as in the proof of \cite[Lemma 11.2]{Naz3}.

\begin{lem}\label{lemma41}
For any integers $k,l\geqslant 0$ and $s_1,\ldots ,s_k, r_1,\ldots ,r_l\geqslant 1$ and for any indices $i_1,\ldots ,i_{k+l},j_1,\ldots ,j_{k+l}=1,\ldots,m+n$ the following implication holds:
$$
\text{if}\quad
\big<t_{i_1 j_1}^{(s_1)}\ldots t_{i_k j_k}^{(s_k)}, t_{i_{k+1} j_{k+1}}^{(-r_1)}\ldots t_{i_{k+l} j_{k+l}}^{(-r_l)}\big>\neq 0,
\quad\text{then}\quad
s_1+\ldots + s_k\geqslant r_1+\ldots +r_l.
$$
\end{lem}

We can now define a bilinear pairing
$
\gr_1 \Y\times\gr_2\Yd\to \CC
$
by 
\beq\label{pairing}
\big<\hat{t}_{i_1 j_1}^{(s_1)}\ldots \hat{t}_{i_k j_k}^{(s_k)}, \bar{t}_{i_{k+1} j_{k+1}}^{(-r_1)}\ldots \bar{t}_{i_{k+l} j_{k+l}}^{(-r_l)}\big>
=
\delta_{s_1+\ldots + s_k, r_1+\ldots +r_l}
\big<t_{i_1 j_1}^{(s_1)}\ldots t_{i_k j_k}^{(s_k)}, t_{i_{k+1} j_{k+1}}^{(-r_1)}\ldots t_{i_{k+l} j_{k+l}}^{(-r_l)}\big>.
\eeq
 For any $r \geqslant 1$ denote by $\gr_{1,r} \Y$ (resp. $\gr_{2,-r}\Yd$) the subspace of degree $r$ (resp. $-r$) in the  graded algebra $\gr_{1 } \Y$ (resp. $\gr_{2 }\Yd$).

\begin{thm}
The map \eqref{mapYdx} is an isomorphism of $\ZZ_2$-graded algebras.
\end{thm}

\begin{prf}
The proof of the theorem can be carried out in parallel with  the case of   dual Yangian for the  general linear Lie algebra $\mathfrak{gl}_N$, as given in  \cite[Sect. 12]{Naz3}. In particular, it uses  the nondegeneracy of restriction of bilinear pairing \eqref{pairing} to $\gr_{1,s} \Y\times \gr_{2,-s}\Yd$   for all $s\geqslant 0$, which can be verified by   similar arguments as \cite[Prop. 12.2]{Naz3}. Roughly speaking, this  is due to the fact that   the $R$-matrices  which govern the defining relations for the  (dual) Yangian  for $\mathfrak{gl}_N$ and for $\mathfrak{gl}_{m|n}$ are of the same form. 
\end{prf}

\subsection{Double Yangian for \texorpdfstring{$\glmn$}{glm|n} at the level $0$}\label{subsec0103}

First, we shall derive a certain lemma on evaluation representations of  the loop Lie superalgebra $\mathcal{L}(\glmn)=\glmn\ot\CC[t,t^{-1}]$, which we need in the proof of Theorem \ref{thmlvlz} below. Recall that for any   representation $\sigma$ of $\glmn$ on $ \CC^{m|n}$ and nonzero   $a\in \CC $ one can define  the evaluation representation    $\sigma_a \colon  \mathcal{L}(\glmn)  \to \ndo\CC^{m|n}$ by
$$
\sigma_a \colon
x\ot t^r \mapsto a^r \sigma (x)\qquad\text{for all } x\in \glmn\text { and }r\in\ZZ.
$$
Denote by $\sigma_{a_1,\ldots ,a_k}$   the tensor product of evaluation representations $\sigma_{a_1},\ldots ,\sigma_{a_k}$. We use the same notation for the  extension of $\sigma_{a_1,\ldots ,a_k}$ to the representation  of    enveloping algebra $\U(\mathcal{L}(\glmn))$. The next lemma  is verified by arguing as in the proof of \cite[Prop 2.2]{Naz2}.

\begin{lem}\label{kernels}
Let $\sigma$ be a faithful representation of $\glmn$. The intersection of all kernels of representations $\sigma_{a_1,\ldots ,a_k}$ with $k> 0$ and nonzero   $a_1,\ldots ,a_k\in \CC$ in  $\U(\mathcal{L}(\glmn))$ is trivial.
\end{lem} 

Throughout the rest of the paper, we assume that $m\neq n$.
The
{\em double Yangian $\dyo$ for $\mathfrak{gl}_{m|n}$ at the level 0} is defined as the $\ZZ_2$-graded unital associative algebra generated by the elements $t_{ij}^{(r)}$ and $t_{ij}^{(-r)}$, where $  i,j=1,\ldots , m+n$ and $r=1,2,\ldots,$ subject to the  defining relations which are written in terms of the generator matrices \eqref{teu} and \eqref{teud}.
They are given by \eqref{YRTT}, \eqref{YdRTT} and
\beq\label{DYRTT}
R(u-v)\ts T_1 (u)\ts T_2^+ (v)
=T_2^+ (v)\ts T_1 (u)\ts R(u-v).
\eeq
The parity of
the elements $t_{ij}^{(r)}$ and
$t_{ij}^{(-r)}$ is again $\bar{i}+\bar{j}$.

The degree operator $\deg_2$, as given by \eqref{deg21} and \eqref{deg22}, defines an ascending filtration 
\beq\label{filtration}
\ldots \subseteq\dyo^{(r)}\subseteq \dyo^{(r+1)}\subseteq\ldots\subseteq\dyo,
\eeq
where $\dyo^{(r)}$ is the linear span of the elements of $\dyo$ whose degrees do not exceed $r$.
 We shall write $\bar{t}_{ij}^{(\pm r)}$  for the images of the generators    in the respective components of the graded  algebra   $\gr_2 \dyo$. 
A direct calculation shows that the   assignments \eqref{assignments1} and  \eqref{assignments2} define a homomorphism
\beq\label{mapYd}
\U(\mathcal{L}(\glmn))\to \textstyle\gr_2  \dyo.
\eeq

Let $ {\rm Y}$ (resp. ${\rm Y}^+$) be the  unital   subalgebra of the double Yangian generated by   all elements $t_{ij}^{(r)}$ (resp. $t_{ij}^{(-r)}$) with $i,j=1,\ldots ,m+n$ and $r=1,2,\ldots .$ Consider the descending filtration on ${\rm Y}^+$ defined by setting the degree of $t_{ij}^{(-r)}$ to be $r$ and denote by $\wtld{{\rm Y}}^+$ the corresponding completion of ${\rm Y}^+$.
Introduce the {\em extended double Yangian  $\dyoc$ at the level 0} as the space of all finite linear combinations of all products  $xy$ for $x\in \wtld{\rm Y}^+$ and $y\in {\rm Y}$ with the multiplication extended by continuity from the double Yangian.
The Hopf superalgebra structure on $\dyoc$ is defined by the formulae
\begin{align}
&\Delta(t_{ij}(u))=\sum_{k=1}^{m+n} t_{ik}(u)\ot t_{kj}(u),\quad
\Delta(t^+_{ij}(u))=\sum_{k=1}^{m+n} t^+_{ik}(u)\ot t^+_{kj}(u),\label{antipodantipod}
\\
&S(T(u))=T(u)^{-1},\quad 
S(T^+(u))=T^+ (u)^{-1}, \qquad
\varepsilon (T(u))= \varepsilon (T^+(u))=1. \label{antipod}
\end{align}

\begin{lem}
For any nonzero $a\in\CC$ the   assignments
$$
t_{ij}^{(r)}\mapsto (-1)^{\bar{i}} a^{r-1} e_{ij}\fand 
t_{ij}^{(-r)}\mapsto (-1)^{\bar{i}} a^{-r} e_{ij}
$$
 with $i,j=1,\ldots, m+n$  and $r=1,2,\ldots$ define a representation of the double Yangian
\beq\label{pia}
\pi_a \colon \dyo \to \ndo \CC^{m|n}.
\eeq
\end{lem}

\begin{prf}
One easily checks that the map \eqref{pia} preserves the ideal of defining relations for the double Yangian  using the equivalent expressions for the above assignments,
$$
T(u)\mapsto R(-u+a)^{\tau_1}\fand
T^+(u)\mapsto R(a-u )^{\tau_1},
$$
where $\tau_1$ is the   transposition $  e_{ij}\mapsto (-1)^{\bar{i}\bar{j}+\bar{i}}e_{ji}$ applied on the first tensor factor.
\end{prf}

Let us define a total order  on the   double Yangian generators as follows.  For any $i,j,k,l=1,\ldots ,m+n$ and $r,s=1,2,\ldots  $   set
\begin{enumerate}[label={(\roman*)}]
\item\label{reef1}   $t_{ij}^{(-r)}\prec t_{kl}^{(s)}$;
\item\label{reef2}  $t_{ij}^{(\pm r)}\prec t_{kl}^{(\pm s)}$ if  $(i,j)$ precedes $(k,l)$ in lexicographical order;
\item\label{reef3} $t_{ij}^{(r)}\prec t_{ij}^{(s)}$ and $t_{ij}^{(-s)}\prec t_{ij}^{(-r)}$ if $r<s$.
\end{enumerate} 
Let us equip the algebra $\dyo$ with the topology     induced by    filtration \eqref{filtration}. 

\begin{thm}\label{thmlvlz}
The  set of all ordered monomials in the generators, with at most power $1$ for odd generators, forms  a topological basis of $\dyo$. 
Hence, the map \eqref{mapYd} is an isomorphism of $\mathbb{Z}_2$-graded algebras.
\end{thm}

\begin{prf}
Clearly, the double Yangian is spanned by all monomials in the generators. Hence, it is sufficient to check that for any integer $p$  and a monomial $\mu$ we can write $\mu $  modulo $\dyo^{(p)}$ as a linear combination of ordered monomials, with at most power $1$ for odd generators. 
First, as with the even case, one  proves by induction using the defining relations that $\mu$ can be expressed as a linear combination of monomials in generators satisfying \ref{reef1} and \ref{reef2}. These monomials are products of submonomials of the form
$t_{ij}^{(\pm r_1)}\ldots t_{ij}^{(\pm r_k)}$. 
If $\bar{i}+\bar{j}$ is even,
  by relations \eqref{Y} and \eqref{Yd}, we have $ [t_{ij}^{(\pm r)},t_{ij}^{(\pm s)} ]=0$, so we can assume that \ref{reef3} holds for such submonomials. If $\bar{i}+\bar{j}$ is odd and $r\neq s$,
   relations \eqref{Y} and \eqref{Yd} take the form
	$$
	 t_{ij}^{(\pm r)} \ts t_{ij}^{(\pm s)}  = -  t_{ij}^{(\pm s)} \ts t_{ij}^{(\pm r)} +\text{lower degree monomials}.
	$$
Hence, if $\pm r_{\sigma (1)}\leqslant\ldots \leqslant \pm r_{\sigma (k)}$ for some permutation $\sigma$ of the indices  $1,\ldots ,k$,
we have
$$
t_{ij}^{(\pm r_1)}\ldots t_{ij}^{(\pm r_k)}
=\varepsilon\ts t_{ij}^{(\pm r_{\sigma(1)})}\ldots t_{ij}^{(\pm r_{\sigma(k)})}
+\text{lower degree monomials}\quad\text{for some }\varepsilon\in\left\{-1,1\right\}.
$$
Finally, we can exclude all   squares of odd generators. Indeed,   by setting    $r=s$ and $(k,l)=(i,j)$ with $\bar{i}+\bar{j}$ odd in relations \eqref{Y} and \eqref{Yd}, on the left hand-side we obtain $2(t_{ij}^{(\pm r)} )^2$, while the right hand-side consists of elements of lower degrees. 
Let us assume that the original monomial $\mu$ is of degree $d>p$. By the preceding discussion, we can express it modulo $\dyo^{(d-1)}$ as a linear combination of ordered monomials of degree    $d$ in the generators, with at most power $1$ for odd generators. Clearly, we can now  continue inductively to express $\mu$ modulo $\dyo^{(p)}$ as a linear combination of such monomials of degrees $d,d-1,\ldots , p+1$, as required.

As for the linear independence, it is verified by using the ideas of Etingof and Kazhdan \cite[Prop. 3.15]{EK3} and Nazarov \cite[Prop. 2.2]{Naz2} which rely on the existence of  evaluation representations. Suppose that some nontrivial linear combination of  ordered monomials   vanishes. Then its image under $\pi_{a_1}\ot\ldots \ot\pi_{a_k}$ is trivial for any choice of $k>0$ and nonzero $a_1,\ldots ,a_k\in\CC$. This leads to a contradiction  by regarding the top degree components of the monomials with respect to the filtration \eqref{filtration} and arguing as in \cite[Prop. 2.2]{Naz2}. In particular, the argument employs the corresponding evaluation representations $\sigma_{a_1,\ldots ,a_k}$ of $\U(\mathcal{L}(\glmn))$ and Lemma \ref{kernels}.
\end{prf}

\begin{rem}
Note that the basis from Theorem \ref{thmlvlz} is topological. For example, the expression for $(t_{ij}^{(-1)})^2$, when $\bar{i}+\bar{j}$ is odd, modulo $\dyo^{(-p)}$, $p\geqslant 4$ is a linear combination  of the monomials  $t_{ij}^{(-r)}t_{ij}^{(-1)}$, $r=2,\ldots ,p-2$ with all  coefficients nonzero.
\end{rem}

\subsection{Double Yangian for  \texorpdfstring{$\mathfrak{gl}_{m|n}$}{glm|n}}\label{subsec0104}

Consider the $R$-matrix \eqref{Rmatrix}. As $P^2=1$, it satisfies  
\beq\label{preuni1}
R(u) R(-u)=1-u^{-2}.
\eeq
Let $g(u) $ be the unique formal power series in  $1+u^{-1}\CC[[u^{-1}]]$ such that
\beq\label{functg}
g(u+m-n)=(1-u^{-2})g(u).
\eeq
The series $g(u)$ also satisfies
\beq\label{preuni2}
g(u) g(-u)(1-u^{-2})=1;
\eeq
see \cite[Sect. 2.2]{JKMY} for more details on  $g(u)$. By combining \eqref{preuni1} and \eqref{preuni2} we find that the normalized $R$-matrix 
$ 
\R(u)=g(u)R(u)
$ 
 possesses the {\em unitarity property},
\beq\label{uni}
\R(u) \ts\R(-u)=1.
\eeq
Also, it is clear that it   satisfies Yang--Baxter equation \eqref{ybe}.
We explain the motivation for the use of this particular normalizing function $g(u)$ in Subsection \ref{subsec00105}; cf. \eqref{rmatpropc}.

The next lemma generalizes \cite[Lemma 2.1]{EK} to the super setting. As with its even counterpart, it is proved by a direct calculation  relying on      $RTT$-relations \eqref{YRTT} and \eqref{YdRTT}.

\begin{lem}\label{plemma21}
For any $c\in\CC$ there exists a unique action of the super Yangian $\Y$ on the algebra $\Yd$ such that for any integer $k\geqslant 0$ we have
\beq\label{lemma21}
T_0 (u) \ts T_1^+(v_1)\ldots T_k^+(v_k)
= (\R_{01}^+)^{-1}\ldots (\R_{0k}^+)^{-1} \ts T_1^+(v_1)\ldots T_k^+(v_k)\ts  \R_{0k}^- \ldots \R_{01}^-
\eeq
on the tensor product
\beq\label{tp}
\ndo\CC^{m|n}\ot
(\ndo\CC^{m|n})^{\ot k}
\ot \Yd
,\eeq 
where the $R$-matrices $\R_{0j}^\pm =\R_{0j}(u-v_j\pm c/2)$ are applied on the tensor factors $0$ and $j$ of \eqref{tp}. In particular, for $k= 0$ we have the identity $T(u)1=1$ on 
$\ndo\CC^{m|n}
\ot \Yd$.
\end{lem}

The
{\em double Yangian $\dyg$ for $\mathfrak{gl}_{m|n}$} is defined as the $\ZZ_2$-graded unital associative algebra generated by the elements $C$, $t_{ij}^{(r)}$ and $t_{ij}^{(-r)}$, where $  i,j=1,\ldots , m+n$ and $r=1,2,\ldots,$ subject to the  defining relations which are written in terms of the generator matrices \eqref{teu} and \eqref{teud}.
They are given by \eqref{YRTT}, \eqref{YdRTT} and 
\beq\label{DYRTTC}
\R(u-v+C/2)\ts T_1 (u)\ts T_2^+ (v)
=T_2^+ (v)\ts T_1 (u)\ts \R(u-v-C/2),
\eeq
where $C$ is   even central element and, as before, the parity of
the elements $t_{ij}^{(\pm r)}$  is  $\bar{i}+\bar{j}$.

There exists a natural epimorphism $\phi\colon \dyg\to \dyo$ such that $t_{ij}^{(\pm r)}\mapsto t_{ij}^{(\pm r)}$ and $C\mapsto 0$.
Due to Theorem \ref{thmlvlz}, the     subalgebra of the double Yangian $\dyg$ generated by all elements $t_{ij}^{(r)}$ (resp. $t_{ij}^{(-r)}$) with $i,j=1,\ldots ,m+n$ and $r=1,2,\ldots  $ coincides with the Yangian  $\Y$ (resp. dual Yangian $\Yd$). 
Consider the descending filtration on the dual Yangian defined by setting the degree of $t_{ij}^{(-r)}$ to be $r$. Denote by $\wtld{{\rm Y}}^+(\glmn)$ the corresponding completion of $\Yd$. We  refer to $\wtld{{\rm Y}}^+(\glmn)$ as the {\em extended dual Yangian}. 
Define the {\em extended double Yangian  $\dyc$} as the space of all finite $\CC[C]$-linear combinations of all products  $xy$ for $x\in\Ydt$ and $y\in \Y$ with the multiplication extended by continuity from the double Yangian. The Hopf superalgebra structure on $\dyc$ is given by
\begin{align*}
&\Delta(t_{ij}(u))=\sum_{k=1}^{m+n} t_{ik}(u+C_2/4)\ot t_{kj}(u-C_1/4),\\
&\Delta(t^+_{ij}(u))=\sum_{k=1}^{m+n} t^+_{ik}(u-C_2/4)\ot t^+_{kj}(u+C_1/4),
\\
&\Delta(C)=C_1 +C_2,\quad
S(C)=-C,\quad 
\varepsilon ( C)= 0 
\end{align*}
and \eqref{antipod}, where $C_1=1\ot C$ and $C_2=1\ot C$.

Extend the degree operator $\deg_2$,   given by \eqref{deg21} and \eqref{deg22}, by defining the degree of the central element $C$ to be zero. Thus, we obtain the ascending filtration 
\beq\label{filtration2}
\ldots \subseteq\dyg^{(r)}\subseteq \dyg^{(r+1)}\subseteq\ldots\subseteq\dyg,
\eeq
where $\dyg^{(r)}$ is the linear span of the elements of $\dyg$ whose degrees do not exceed $r$.
  Denote the corresponding graded algebra by $\gr_2 \dyg$ and denote the images of   generators in its respective components by $\bar{t}_{ij}^{(\pm r)}$ and $\bar{C}$. A direct calculation  relying on the supercommutation relations \eqref{glmn}  
and the
form  of the series $g(u)$, which goes in parallel with the proof of \cite[Prop. 2.1]{JKMY},  implies that 
 the assignments
$$
e_{ij}(r-1)\mapsto (-1)^{\bar{i}}\ts \bar{t}_{ij}^{(r)},\qquad e_{ij}(-r)\mapsto (-1)^{\bar{i}}\ts \bar{t}_{ij}^{(-r)}   \Fand K\mapsto \bar{C}
$$
with
$i,j=1,\ldots ,m+n$  and $r\geqslant 1$ define a homomorphism
\beq\label{mapYYd}
\U(\glmnht) \to \textstyle\gr_2 \dyg.
\eeq

 Extend the total order   \ref{reef1}--\ref{reef3}  to $ \dyg$ so that it includes the central element $C$ in an arbitrary way. In the following theorem, we consider the topology over $\dyg$   induced by the filtration \eqref{filtration2}.

\begin{thm}\label{pbbwmain}
The  set of all ordered monomials in the generators, with at most power $1$ for odd generators, forms  a topological basis of $\dyg$. 
Hence, the map \eqref{mapYYd} is an isomorphism of $\mathbb{Z}_2$-graded algebras.
\end{thm}

\begin{prf}
As with the level $0$ case from Theorem \ref{thmlvlz},
 one     verifies by induction  using   the defining relations  that  the ordered monomials in   generators, with at most power $1$ for odd generators, topologically span $ \dyg$.  
Next, Lemma \ref{plemma21} implies that, for any $c\in\CC$, the dual Yangian is naturally equipped with the structure of module for the double Yangian of level $c$, i.e. such that $C$ acts as a scalar multiplication by $c$. Thus, the central element $C$ is nonzero. Moreover,   by arguing as in the proof of \cite[Thm. 2.2]{JKMY}, one can show that its powers $1, C,\ldots ,C^k$ are linearly independent for any positive $k$. 
Finally, the linear independence of the ordered monomials in the generators with at most power $1$ for odd generators is established by   examining the image of their linear combination under the map $(1\ot\phi)\circ\Delta$ and using Theorem \ref{thmlvlz},  in parallel with the
 corresponding part of  the proof of \cite[Thm. 2.2]{JKMY}. 
\end{prf}

\subsection{Quantum contraction}\label{subsec00105}

In the following lemma, we introduce  the  dual Yangian counterpart of the series $z (u)$; recall \eqref{qcontra}. 

\begin{lem}\label{wldfndlm}
There exists a unique power series $z^+ (u)$ in
$\wtld{{\rm Y}}^+(\glmn)[[u]]$ such that 
\beq\label{qcontrad}
P^{\tau_2}\ts  T^+_1(u+m-n) \left(T^+_2(u)^{-1}\right)^\tau = P^{\tau_2}\ot  z^+ (u).
\eeq
\end{lem}

\begin{prf}
First of all, we observe that  the coefficients of  matrix entries of the shifted series $T^+(u+m-n)$ and the inverse $T^+(u)^{-1}$ are well-defined elements of the extended dual Yangian $\wtld{{\rm Y}}^+(\glmn)$.
The lemma is proved by analogous arguments as its Yangian counterpart established in  \cite[Sect. 1]{Naz}, but we provide some details for completeness. First, consider the defining relation \eqref{YdRTT}. Multiplying it from the right and from the left by $T_2^+(v)^{-1}$ and then applying the transposition $\tau$ on the second tensor component we obtain
$$
R(u-v)^{\tau_2}  \left(T_2^+ (v)^{-1}\right)^\tau T_1^+ (u)
=  T_1^+ (u) \left(T_2^+ (v)^{-1}\right)^\tau R(u-v)^{\tau_2}.
$$
Clearly, this is equivalent with
\beq\label{eeq1}
 \left(R(u-v)^{\tau_2}\right)^{-1} T_1^+ (u) \left(T_2^+ (v)^{-1}\right)^\tau 
=  \left(T_2^+ (v)^{-1}\right)^\tau T_1^+ (u)\left(R(u-v)^{\tau_2}\right)^{-1}.
\eeq
As $(P^{\tau_2})^2=(m-n)P^{\tau_2}$, we have
$$
\left(R(u)^{\tau_2}\right)^{-1}=\sum_{l\geqslant 0}\frac{\left(P^{\tau_2}\right)^l}{u^l}
=1+\sum_{l\geqslant 1}\frac{\left(m-n\right)^{l-1}}{u^l}P^{\tau_2}
=1+(u-m+n)^{-1}P^{\tau_2}.
$$
Hence, multiplying \eqref{eeq1} by $u-v-m+n$ we get
 $$
\left(u-v-m+n +P^{\tau_2}\right)  T_1^+ (u) \left(T_2^+ (v)^{-1}\right)^\tau 
=  \left(T_2^+ (v)^{-1}\right)^\tau T_1^+ (u)\left(u-v-m+n +P^{\tau_2}\right).
$$
Replacing the variables $(u,v)$ by $(u+m-n,u)$ the above identity becomes
\beq\label{eeq2}
 P^{\tau_2}\ts   T_1^+ (u+m-n) \left(T_2^+ (u )^{-1}\right)^\tau 
=  \left(T_2^+ (u )^{-1}\right)^\tau T_1^+ (u+m-n)\ts  P^{\tau_2}.
\eeq
Finally, we observe that the image  of   
$ 
P^{\tau_2}=\sum_{i,j=1}^{m+n} e_{ij}\ot e_{ij}\ts (-1)^{\bar{i}\bar{j}}
$ 
is one-dimensional, so that  the equality \eqref{eeq2} implies the assertion of the lemma. 
\end{prf}

As with the original quantum contraction, the series $z^+(u) $ gives rise to a family of central elements in the extended dual Yangian $\wtld{{\rm Y}}^+(\glmn)$.

\begin{lem}\label{josjednalema}
All coefficients of $z^+ (u)$ belong to the centre of the   algebra $\wtld{{\rm Y}}^+(\glmn)$.
\end{lem}

\begin{prf}
The   lemma can be verified by suitably modifying   the proof of \cite[Thm. 1]{Naz}. 
In the proof of Theorem \ref{propqconttra} below we already present in detail such arguments  in the case of the normalized $R$-matrix, so here we     only  sketch  the major steps of the proof. 
First, by   \eqref{YdRTT} and \eqref{eeq1} we have
\begin{align}
&P_{12}^{\tau_2}
\left(R_{02}(v-u )^{\tau_2}\right)^{-1}
R_{01}(v-u-m+n)\ts
T^+_0(v)\ts T_1^+(u+m-n)
\left(T_2^+(u )^{-1}\right)^\tau\non\\
&\quad=
P_{12}^{\tau_2}\ts
T_1^+(u+m-n)
\left(T_2^+(u )^{-1}\right)^\tau
T^+_0(v)\ts
\left(R_{02}(v-u )^{\tau_2}\right)^{-1}
R_{01}(v-u-m+n).\label{alfnjfd8}
\end{align}
Next, using the the property
\beq\label{rmatprop}
P_{12}^{\tau_2}
\left(R_{02}(u )^{\tau_2}\right)^{-1}
R_{01}(u-m+n)
=P_{12}^{\tau_2} \left(1- (u-m+n)^{-2}\right) 
\eeq
of the $R$-matrix \eqref{Rmatrix} and Lemma \ref{wldfndlm}, we bring \eqref{alfnjfd8}  to the form
$$
\left(1-(v-u-m+n)^{2}\right) P_{12}^{\tau_2} \ts
T^+_0(v)\ts z^+ (u)
=
\left(1-(v-u-m+n)^{2}\right) P_{12}^{\tau_2} \ts
 z^+ (u)
T^+_0(v) .
$$
Finally, we cancel  the terms $ 1-(v-u-m+n)^{2} $ and conclude  that $T^+_0(v)$ and  $z^+ (u)$ commute, as required.
\end{prf}

Observe that the proof of Lemma \ref{josjednalema} uses the property
\eqref{rmatprop} of  $R(u)$. By combining \eqref{functg} and \eqref{rmatprop}, one finds that the normalized $R$-matrix $\R(u)$ satisfies the identity
\beq\label{rmatpropc}
P_{12}^{\tau_2}
\left(\R_{02}(u )^{\tau_2}\right)^{-1}
\R_{01}(u-m+n)
=P_{12}^{\tau_2}, 
\eeq
   a super analogue of the ordinary {\em crossing symmetry property} for the even Yang $R$-matrix. 
This property plays a key role in the proof of the next theorem.

\begin{thm}\label{propqconttra}
All coefficients of the series $z(u)$ and $z^+(u)$ belong to the centre of the extended double Yangian $\dyc$.
\end{thm}

\begin{prf}
Clearly, to prove the theorem, it is sufficient to verify the equalities
\beq\label{dydvaide}
T (v)\ts z^+(u) = z^+(u)\ts T (v)
\fand
T^+(v)\ts z(u) = z(u)\ts T^+(v).
\eeq
Let us prove the first equality. Set $x=v-u+C/2$. Consider the identity
\begin{align}
&P_{12}^{\tau_2}
\left(\R_{02}(x)^{\tau_2}\right)^{-1}
\R_{01}(x-m+n )\ts
T_0(v)\ts T_1^+(u+m-n)
\left(T_2^+(u)^{-1}\right)^\tau\non\\
=&\ts 
P_{12}^{\tau_2}\ts
T_1^+(u+m-n)
\left(T_2^+(u)^{-1}\right)^\tau
T_0(v)
\left(\R_{02}(x -C )^{\tau_2}\right)^{-1}
\R_{01}(x-m+n-C ), \label{algnjfd8}
\end{align}
which is deduced from \eqref{DYRTTC}. 
By using \eqref{rmatpropc} and then  Lemma \ref{wldfndlm}, we rewrite the left-hand side of \eqref{algnjfd8} as
\begin{align}
&P_{12}^{\tau_2}\ts 
T_0(v)\ts T_1^+(u+m-n)
\left(T_2^+(u)^{-1}\right)^\tau 
= 
T_0(v)\ts P_{12}^{\tau_2}\ts T_1^+(u+m-n)
\left(T^+_2(u)^{-1}\right)^\tau\non \\
=   &\ts
T_0(v)\ts P_{12}^{\tau_2}\ts z^+(u)=P_{12}^{\tau_2}\ts T_0(v)\ts  z^+(u).\label{idd1}
\end{align}
As for the right-hand side of \eqref{algnjfd8}, using Lemma \ref{wldfndlm} and then \eqref{rmatpropc} we get
\begin{align}
&
P_{12}^{\tau_2} 
\ts z^+(u)\ts
T_0(v) 
\left(\R_{02}(x -C )^{\tau_2}\right)^{-1}
\R_{01}(x-m+n-C ) \non\\
=&\ts  z^+(u)\ts
T_0(v)\ts P_{12}^{\tau_2} 
\left(\R_{02}(x -C )^{\tau_2}\right)^{-1}
\R_{01}(x-m+n-C )\non\\
=& z^+(u)\ts
T_0(v)\ts P_{12}^{\tau_2} 
= P_{12}^{\tau_2} \ts  z^+(u)\ts
T_0(v).\label{idd2}
\end{align}
Due to \eqref{algnjfd8}, the expressions in \eqref{idd1} and \eqref{idd2} coincide, which implies the first equality in \eqref{dydvaide}. Regarding the second equality, it is verified by an analogous argument,   relying on Lemma \ref{wldfndlm} and   \eqref{rmatpropc}. However, instead of \eqref{algnjfd8}, its proof starts with the identity
\begin{align*}
&P_{12}^{\tau_2}
\left(\R_{02}(y)^{\tau_2}\right)^{-1}
\R_{01}(y-m+n )\ts
T^+_0(v)\ts T_1 (u+m-n)
\left(T_2 (u)^{-1}\right)^\tau \\
=&\ts 
P_{12}^{\tau_2}\ts
T_1 (u+m-n)
\left(T_2 (u)^{-1}\right)^\tau
T^+_0(v)
\left(\R_{02}(y +C )^{\tau_2}\right)^{-1}
\R_{01}(y-m+n+C )   
\end{align*}
with $y=-u+v-C/2$, which again follows from the defining relation \eqref{DYRTTC}.
\end{prf}

We now give a simple application of Theorem \ref{propqconttra}. Introduce the {\em extended vacuum module $\Vccgl$ at the level $c\in\CC$} as the quotient of the algebra  $\dyc$ by its left ideal generated by $C-c\cdot 1$ and the elements $t_{ij}^{(r)}$, where $i,j=1,\ldots ,m+n$ and $r=1,2,\ldots .$   Denote by $\vac$ the image of the unit $1\in \dyc$ in the extended vacuum module. Let
$$
\mathfrak{z} (\Vccgl )
=\left\{
v\in \Vccgl\,:\, t_{ij}^{(r)}v=0\text{ for all }i,j=1,\ldots ,m+n,\, r=1,2,\ldots 
\right\}
$$ 
be the {\em subspace of invariants of $\Vccgl$}.
Theorem \ref{propqconttra} implies
\begin{kor}
All coefficients of the series $z^+ (u)\vac$ belong to the subspace of invariants $\mathfrak{z} (\Vccgl )$ of the extended vacuum module.
\end{kor}

At the end, we point out some similarities and differences between the series
$z (u)$ and $z^+(u)$. First of all, we remark that Lemmas \ref{wldfndlm}, \ref{josjednalema} and the following discussion hold for $m=n$ as well. As we recall  in Section \ref{subsec0101}, the coefficients of $z (u)$ generate the centre of the Yangian. However, its dual Yangian counterpart does not   exhibit such a property. To see this,
write $t^+(u)=1-T^+(u)$. Using the formal Taylor Theorem we find
$$
T_1^+(u+m-n)\left(T_2^+(u)^{-1}\right)^{\tau}
=\left(
1-\sum_{r\geqslant 0} \frac{(m-n)^r}{r!}\frac{d^r}{du^{r}}t_1^+(u)
\right)
\left(
 \sum_{r\geqslant 0}  t_2^+(u)^r
\right)^\tau .
$$
Note that for all $r\geqslant 1$ the degree of the coefficient of $u^{r-1}$ in $\frac{d^r}{du^{r}}t_1^+(u)$ and in $t_2^+(u)^r$ is less than or equal to $-r-1$. Furthermore, the degree of  the coefficient of $u^{r-1}$ in $t_1^+(u)t_2^+(u)^\tau$ is $-r-2$. Thus, the above expression is of the form
$$
T_1^+(u+m-n)\left(T_2^+(u)^{-1}\right)^{\tau}
=
1-t_1^+(u)+t_2^+(u)^\tau+\text{lower degree terms}.
$$
As $P^{\tau_2}t_1^+(u)=P^{\tau_2}t_2^+(u)^\tau$, this implies that the degree of the coefficient $z^{-r}$ in  
$$
z^+(u)=1-\sum_{r\geqslant 1}z^{(-r)} \ts u^{r-1}
$$
does not exceed $-r-1$. Hence, in particular, the coefficients of the series $z^+(u)$ do not contain any central elements  of   degree $-1$. On the other hand,   they are algebraically independent. Indeed, 
extend the degree function $\deg_2$ for the dual Yangian  to the extended algebra $\wtld{{\rm Y}}^+(\glmn)$   by allowing it to take the
infinite value.  Then the elements of finite
degree form a subalgebra which we denote by $\wtld{{\rm Y}}^+(\glmn)_{\text{fin}}$.
All coefficients $z^{(-r)}$ belong to $\wtld{{\rm Y}}^+(\glmn)_{\text{fin}}$ and their  images $\bar{z}^{(-r)}$    in $\gr_2 \wtld{{\rm Y}}^+(\glmn)_{\text{fin}}\cong  \U(t^{-1}\glmn[t^{-1}])$ are given by
$$
\bar{z}^{(-r)}=r\sum_{i=1}^{m+n} e_{ii}\ot t^{-r-1}\quad
\text{for all }r\geqslant 1.
$$

\section{Reflection algebras for  the Lie superalgebra $\mathfrak{gl}_{m|n}$}\label{section2}
In this section, we  study reflection algebras, a certain class of   left coideal subalgebras of the extended double Yangian at the level $0$. In particular, we give a presentation of such algebras   by generators and relations.

\subsection{Reflection subalgebras of  the double Yangian}\label{subsec0105}

For any $\ell=1,\ldots ,m+n$ let $G=(g_{ij})$ be the diagonal matrix of order $m+n$,
\beq\label{gmatricca}
G=\diag(\varepsilon_1,\ldots, \varepsilon_{m+n}), \qquad\text{where}\qquad
\varepsilon_i=\begin{cases}
1&\text{for }i=1,\ldots ,\ell,\\
-1&\text{for }i=\ell+1,\ldots,m+n.
\end{cases}
\eeq
Introduce the matrices
\beq\label{bbplus}
B(u)=T(u )\ts G\ts T(-u)^{-1}\Fand B^+(u)=T^+(u)\ts G\ts T^+(-u)^{-1} 
\eeq
and denote by $b_{ij}(u)$ and $b^+_{ij}(u)$ their matrix entries,
\beq\label{serieess}
B(u)=\sum_{i,j=1}^{m+n} (-1)^{\bar{i}\bar{j}+\bar{j}} e_{ij}\ot b_{ij}(u)
\Fand
B^+(u)=\sum_{i,j=1}^{m+n} (-1)^{\bar{i}\bar{j}+\bar{j}} e_{ij}\ot b^+_{ij}(u).
\eeq
Finally, let $b_{ij}^{(r)}$ and  $b_{ij}^{(-r)}$ be the coefficients of the series $b_{ij}(u)\in\Y[[u^{-1}]]$
and $b^+_{ij}(u)\in \Ydt[[u]]$, so that we have
\beq\label{seriees}
b_{ij}(u)=g_{ij}+\sum_{r\geqslant 1} b_{ij}^{(r)} u^{-r}
\Fand
b^+_{ij}(u)=g_{ij}-\sum_{r\geqslant 1} b_{ij}^{(-r)} u^{r-1}.
\eeq
From now on, we   regard $b_{ij}^{(r)}$ and  $b_{ij}^{(-r)}$ as  elements of the extended double Yangian at the level $0$.
Using the $RTT$-relations \eqref{YRTT}, \eqref{YdRTT} and \eqref{DYRTT}, along with the identity
$$R(u)G_1 R(v)G_2=G_2 R(v) G_1 R(u),$$
one can show  that the matrices  \eqref{bbplus} satisfy the {\em reflection relations}
\begin{align}
R(u-v)\ts B_1(u)\ts R(u+v )\ts B_2(v)
&=B_2(v)\ts R(u+v )\ts  B_1(u)\ts R(u-v),\label{refl1}\\
R(u-v)\ts B^+_1(u)\ts R(u+v)\ts B^+_2(v)
&=B^+_2(v)\ts R(u+v)\ts  B^+_1(u)\ts R(u-v),\label{refl2}\\
R(u-v )\ts B_1(u)\ts R(u+v )\ts B^+_2(v) 
 &=B^+_2(v)\ts R(u+v )\ts  B_1(u)\ts R(u-v ).\label{refl3} 
\end{align}
Furthermore, they possess the {\em unitarity properties}
\beq\label{unidentities}
B(u)\ts B(-u)=1
\Fand
B^+(u)\ts B^+(-u)=1.
\eeq  

Extend the degree function $\deg_2$ for the double Yangian, as defined by \eqref{deg21},  \eqref{deg22} and $\deg_2 C=0$, to the algebra $\dyoc$   by allowing it to take the
infinite value.  Then the elements of finite
degree form a subalgebra which we denote by $\dyoc_{\text{fin}}$.
One easily checks that   $b_{ij}^{(r)}$ and  $b_{ij}^{(-r)}$ belong to $\dyoc_{\text{fin}}$.
Let $\DB$ be its unital subalgebra   generated by  all elements $b_{ij}^{(r)}$ and  $b_{ij}^{(-r)}$. Moreover, let $\B$ (resp. $\Bd$) be its unital subalgebra generated  by all $b_{ij}^{(r)}$   (resp.  $b_{ij}^{(-r)}$) with $r=1,2,\ldots.$ We shall refer to  all these subalgebras as {\em reflection algebras}.

Consider the families $I_0,I_1\subset \left\{1,\ldots ,m+n\right\}^{\times 2}$ of pairs of indices 
\begin{align*}
&I_0=\left\{(i,j)\,:\,1\leqslant i,j\leqslant\ell\text{ or }\ell +1\leqslant i,j\leqslant m+n\right\},\\
&I_1=\left\{(i,j)\,:\,1\leqslant i \leqslant\ell <j\leqslant m+n\text{ or }1\leqslant j \leqslant\ell <i\leqslant m+n\right\}.
\end{align*}
 Define the families of elements
\begin{align*}
&\Gamma =\left\{b_{ij}^{(2r-1)}\,:\, (i,j)\in I_0\text{ and }r\geqslant 1\right\} \cup 
\left\{b_{ij}^{( 2r )}\,:\, (i,j)\in I_1\text{ and }r\geqslant 1\right\} \subset \B,\\
&\Gamma^+=\left\{b_{ij}^{(-2r)}\,:\, (i,j)\in I_0\text{ and }r\geqslant 1\right\}\cup
\left\{b_{ij}^{(-2r+1)}\,:\, (i,j)\in I_1\text{ and }r\geqslant 1\right\} \subset \Bd.
\end{align*}
By a direct computation we obtain from \eqref{bbplus} the equalities 
\beq\label{slikeizo}
b_{ij}^{(r)}=\left((-1)^{r-1}\varepsilon_i +\varepsilon_j\right)t_{ij}^{(r)}+\ldots
\Fand
b_{ij}^{(-r)}=((-1)^r\varepsilon_i +\varepsilon_j)t_{ij}^{(-r)}+\ldots 
\eeq
for all $i,j=1,\ldots ,m+n$ and $r=1,2,\ldots ,$ 
where the ellipses stand for the lower degree terms, i.e the terms of degree strictly less than $r-1$ (resp. $-r$), of parity   $\bar{i}+\bar{j}$. Hence, in particular, all elements $b_{ij}^{(\pm r)}$ are homogeneous.
Moreover, for all $b_{ij}^{(\pm r)}\in \Gamma\cup \Gamma^+$ we have
\beq\label{stupnjevvi}
\deg_2\ts b_{ij}^{(r)}=r-1\fand
\deg_2\ts b_{ij}^{(-r)}=-r,
\eeq
while  for the remaining elements $b_{ij}^{(\pm r)}\not\in \Gamma\cup \Gamma^+$  we  have 
\beq\label{stupnjevvi2}
\deg_2\ts b_{ij}^{(r)}<r-1\fand
\deg_2\ts b_{ij}^{(-r)}<-r.
\eeq
Introduce the following subsets of $    \gr_2 \dyoc_{\text{fin}} =\gr_2 \dyo \cong \U(\mathcal{L}(\glmn))$:
$$
\wvr{\Gamma}=\left\{\bar{b}_{ij}^{(r)}\,:\, b_{ij}^{(r)}\in \Gamma\right\} 
\fand
\wvr{\Gamma}^+=\left\{\bar{b}_{ij}^{(-r)}\,:\, b_{ij}^{(-r)}\in \Gamma^+\right\}.
$$

\begin{thm}\label{reflprro}
The family 
$\wvr{\Gamma} $ (resp. $\wvr{\Gamma}^+$) generates   $\gr_2 \B$ (resp. $\gr_2 \Bd$). Hence, the union  $\wvr{\Gamma}\cup \wvr{\Gamma}^+ $ generates the algebra $\gr_2 \DB$.
\end{thm}

\begin{prf}
Let us prove that the family  $\wvr{\Gamma}^+$ generates $\gr_2  \Bd$. 
By extracting the constant term of the second identity in \eqref{unidentities} we find
\beq\label{lkmbghdf}
\left(\varepsilon_i +\varepsilon_j\right) b_{ij}^{(-1)}
=\sum_{a=1}^{\ell}b_{ia}^{(-1)}b_{aj}^{(-1)} + \sum_{a=\ell+1}^{m+n}b_{ia}^{(-1)}b_{aj}^{(-1)}.
\eeq
Therefore, if $b_{ij}^{(-1)}$ does not belong to $\Gamma^+$, we see by \eqref{stupnjevvi} and \eqref{stupnjevvi2} that the degree of the right-hand side of \eqref{lkmbghdf} is $-2$. Moreover, we have the equality
$$
\bar{b}_{ij}^{(-1)}
=\begin{cases}
  \frac{1}{2}\sum_{a=\ell+1}^{m+n}\bar{b}_{ia}^{(-1)}\bar{b}_{aj}^{(-1)},&\text{if }1\leqslant i,j\leqslant \ell,\\
-\frac{1}{2}\sum_{a=1}^{\ell}\bar{b}_{ia}^{(-1)}\bar{b}_{aj}^{(-1)},&\text{if }\ell<i,j\leqslant m+n
	\end{cases}
$$
in the corresponding graded algebra $\gr_2  \Bd$. In other words, we proved that all elements $\bar{b}_{ij}^{(-1)}$ belong to the   subalgebra generated by $\wvr{\Gamma}^+$.

By extracting the coefficient of $u$ in the second identity in \eqref{unidentities} we find
\beq\label{lkmbghdf2}
\left(\varepsilon_i -\varepsilon_j\right) b_{ij}^{(-2)}
=
\sum_{a=1}^{\ell}b_{ia}^{(-1)}b_{aj}^{(-2)}
 + \sum_{a=\ell+1}^{m+n}b_{ia}^{(-1)}b_{aj}^{(-2)}
-\sum_{a=1}^{\ell}b_{ia}^{(-2)}b_{aj}^{(-1)} 
- \sum_{a=\ell+1}^{m+n}b_{ia}^{(-2)}b_{aj}^{(-1)}.
\eeq
Therefore, if $b_{ij}^{(-2)}$ does not belong to $\Gamma^+$, we see by \eqref{stupnjevvi} and \eqref{stupnjevvi2} that the degree of the right-hand side of \eqref{lkmbghdf2} is $-3$. Moreover, we have the equality
$$
\bar{b}_{ij}^{(-2)}
=\begin{cases} 
	-\frac{1}{2}\sum_{a=1}^{\ell}\bar{b}_{ia}^{(-2)}\bar{b}_{aj}^{(-1)}
	+\frac{1}{2} \sum_{a=\ell+1}^{m+n}\bar{b}_{ia}^{(-1)}\bar{b}_{aj}^{(-2)},&\text{if }1\leqslant i \leqslant\ell <j\leqslant m+n,\\
\frac{1}{2}\sum_{a=\ell+1}^{m+n}\bar{b}_{ia}^{(-2)}\bar{b}_{aj}^{(-1)}
-\frac{1}{2}\sum_{a=1}^{\ell}\bar{b}_{ia}^{(-1)}\bar{b}_{aj}^{(-2)},&\text{if }1\leqslant j \leqslant\ell <i\leqslant m+n
	\end{cases}
$$
in $\gr_2 \Bd$. Thus,  all elements $\bar{b}_{ij}^{(-2)}$ belong to the   subalgebra generated by $\wvr{\Gamma}^+$.

Suppose that the elements $\bar{b}_{ij}^{(-r)}$, where $i,j=1,\ldots ,m+n$ and $r=1,\ldots,s$, belong to the   subalgebra generated by $\wvr{\Gamma}^+$. Then, by arguing as above, one can show that all elements $\bar{b}_{ij}^{(-s-1)}$ belong to the   subalgebra generated by $\wvr{\Gamma}^+$ as well. Hence, we conclude by induction that the $\wvr{\Gamma}^+$ generates the entire algebra $\gr_2  \Bd$.

The fact that the family $\wvr{\Gamma} $ generates $\gr_2  \B$ can be verified by suitably modifying the above arguments and employing the first identity in \eqref{unidentities}. In fact,  the constant term of $B(u)$ is 1, which slightly simplifies the proof in this case, so we omit it.  The analogous statement for $\gr_1   {\rm B}(\mathfrak{gl}_{m|0})$ was established by Molev and Ragoucy; see \cite[Sect. 3.1]{MR}.  
\end{prf}

Consider the involutive automorphism $\sigma$ of $\mathfrak{gl}_{m|n}  $ given by
\begin{align*}
\sigma\colon e_{ij} \mapsto \varepsilon_i\ts \varepsilon_j \ts e_{ij},\quad\text{where } i,j=1,\ldots ,m+n.
\end{align*}
It induces the decomposition $\mathfrak{gl}_{m|n}=\mathfrak{gl}_{m|n} (-1) \oplus\mathfrak{gl}_{m|n} (1) $, where $\mathfrak{gl}_{m|n} (\pm 1) $ denotes the eigenspace of $\sigma$ corresponding to the eigenvalue $\pm 1$. Let $\mathfrak{gl}_{m|n}[t,t^{-1}]^{\sigma} \subset \glmnht$ be the Lie  superalgebra of all Laurent   polynomials of the form
$$
 \sum_{i=-p}^q a_i \ot t^i  ,\quad\text{where }p,q\in\ZZ_{\geqslant 0},\,a_{i}\in \mathfrak{gl}_{m|n}((-1)^i) .
$$
By using the
 isomorphism \eqref{mapYYd}  and the identities \eqref{slikeizo} we find
$$
\bar{b}_{ij}^{(r)}=(-1)^{\bar{i}}\left((-1)^{r-1}\varepsilon_i +\varepsilon_j\right)e_{ij}(r-1)
\fand
\bar{b}_{ij}^{(-r)}=(-1)^{\bar{i}}((-1)^r\varepsilon_i +\varepsilon_j)e_{ij} (-r)  
$$
for all $b_{ij}^{(\pm r)}\in \Gamma\cup \Gamma^+$. Therefore, by combining Theorems \ref{pbbwmain} and  \ref{reflprro} we obtain 

\begin{kor}\label{koorolar}  There exists an algebra isomorphism
\beq\label{izomorfizam}
\textstyle\gr_2 \DB\cong \U(\mathfrak{gl}_{m|n}[t,t^{-1}]^{\sigma} ).
\eeq
\end{kor}

\begin{rem}
Corollary  \ref{koorolar} can be directly translated into the $\CC[[h]]$-module setting from  paper \cite{K} by the formal rescaling of the generators and spectral parameter. A direct calculation shows that the images of the odd coefficients  of the double Sklyanin determinant $\wtld{\mathbb{A}}_0(u)$ from \cite[Prop. 3.8]{K}
in the corresponding graded algebra are mapped by the isomorphism    \eqref{izomorfizam} to scalar multiples of $I\ot t^{2r}$ with $r\in\ZZ$. Thus, as with the odd coefficients of the ordinary Sklyanin determinant \cite[Thm. 3.4]{MR}, they form an algebraically independent family.
\end{rem}

The next proposition is a generalization of \cite[Prop. 3.3]{MR}, which states that,  in the even case, i.e. for $n=0$, the subalgebra of the Yangian $ {\rm Y}(\mathfrak{gl}_{m })={\rm Y}(\mathfrak{gl}_{m|0})$  generated by all $b_{ij}^{(r)}$,  where $i,j=1,\ldots ,m$ and $r=1,2,\ldots ,$ is a left coideal in  $ {\rm Y}(\mathfrak{gl}_{m})$.

\begin{pro}\label{coidpropossitio}
The subalgebra $\DB$ is a left coideal in $\dyoc$, i.e. we have
$$
\Delta\left(\DB\right)\subseteq \dyoc\ot\DB.
$$
\end{pro}

\begin{prf}
It is sufficient to check that the images of the generators $b_{ij}^{(\pm r)}$ under the coproduct $\Delta$ belong to $\dyoc\ot\DB $. As with the aforementioned result of Molev and Ragoucy, this is proved by a direct calculation which relies on the explicit formulae \eqref{antipodantipod} for the coproduct. First, as the coproduct  maps $T^+(u)T^+(u)^{-1}$ to $1$, we derive from the second formula in \eqref{antipodantipod} that for all $i,j$ we have
$$
\Delta (t'^+_{ij}(u))=\sum_{k=1}^{m+n}t'^+_{kj}(u)\ot t'^+_{ik}(u),\qquad\text{where}\qquad
T^+(u)^{-1}=\sum_{i,j=1}^{m+n} (-1)^{\bar{i}\bar{j}+\bar{j}}e_{ij}\ot t'^+_{ij}(u).
$$
Hence, using the second formula in \eqref{bbplus}, we get
\begin{align*}
\Delta(b_{ij}^+(u))
&=\Delta\left(\sum_{k=1}^{m+n}\varepsilon_k \ts t_{ik}^+(u)\ts t'^+_{kj}(-u)\right)
= \sum_{k,r,s=1}^{m+n}\varepsilon_k\ts  t_{ir}^+(u)\ts t'^+_{sj}(u)\ot t_{rk }^+(u)\ts t'^+_{ks}(-u) \\
&= \sum_{ r,s=1}^{m+n} t_{ir}^+(u)\ts t'^+_{sj}(u)\ot \sum_{ k=1}^{m+n}t_{rk }^+(u)\ts \varepsilon_k\ts  t'^+_{ks}(-u) 
=\sum_{ r,s=1}^{m+n} t_{ir}^+(u)\ts t'^+_{sj}(u)\ot b^+_{rs}( u),
\end{align*}
which belongs to $\dyoc\ot\DB  [[u]]$, as required. The statement for the series $b_{ij}(u)$ is verified analogously; cf. \cite[Prop. 3.3]{MR}.
\end{prf}

The proof of Proposition \ref{coidpropossitio} implies
\begin{kor}
The subalgebra $\Bd$ is a left coideal in the algebra $\wtld{{\rm Y}}^+(\glmn)$.
\end{kor}

\subsection{Presentation of reflection algebras}\label{subsec0106}

Let $G=(g_{ij})$ be the diagonal matrix as in \eqref{gmatricca}.
Define $\DBB$ as the $\ZZ_2$-graded  unital associative algebra generated by the elements   $\beta_{ij}^{(r)}$ and $\beta_{ij}^{(-r)}$, where $i,j=1,\ldots ,m+n$ and $r=1,2,\ldots .$   The parity of   $\beta_{ij}^{(\pm r)}$ is $\bar{i}+\bar{j}$ and the generators are  subject to the   defining relations which are given as follows.
First, we introduce the   series
$$
\beta_{ij}(u)=g_{ij}+\sum_{r\geqslant 1} \beta_{ij}^{(r)} u^{-r}
\Fand
\beta^+_{ij}(u)=g_{ij}-\sum_{r\geqslant 1} \beta_{ij}^{(-r)} u^{r-1} 
$$
with $i,j=1,\ldots ,m+n$. Next, we organize the series  into matrices
$$
\mathcal{B}(u)=\sum_{i,j=1}^{m+n} (-1)^{\bar{i}\bar{j}+\bar{j}} e_{ij}\ot \beta_{ij}(u)
\Fand
\mathcal{B}^+(u)=\sum_{i,j=1}^{m+n} (-1)^{\bar{i}\bar{j}+\bar{j}} e_{ij}\ot \beta^+_{ij}(u).
$$
Finally, the defining relations consist of three {\em reflection equations},
\begin{align}
R(u-v)\ts \mathcal{B}_1(u)\ts R(u+v )\ts \mathcal{B}_2(v)
&=\mathcal{B}_2(v)\ts R(u+v )\ts  \mathcal{B}_1(u)\ts R(u-v),\label{brefl1}\\
R(u-v)\ts \mathcal{B}^+_1(u)\ts R(u+v)\ts \mathcal{B}^+_2(v)
&=\mathcal{B}^+_2(v)\ts R(u+v)\ts  \mathcal{B}^+_1(u)\ts R(u-v),\label{brefl2}\\
 R(u-v )\ts \mathcal{B}_1(u)\ts  R(u+v )\ts \mathcal{B}^+_2(v) 
&=\mathcal{B}^+_2(v)\ts  R(u+v )\ts  \mathcal{B}_1(u)\ts  R(u-v ) \label{brefl3} 
\end{align}
and two {\em unitarity conditions},
\beq\label{bunidentities}
\mathcal{B}(u)\ts \mathcal{B}(-u )=1
\Fand
\mathcal{B}^+(u)\ts \mathcal{B}^+(-u)=1.
\eeq

Denote by $\BB$ (resp. $\BBd$) a unital subalgebra of $\DBB$ generated by   all $\beta_{ij}^{(r)}$   (resp.  $\beta_{ij}^{(-r)}$) with $i,j=1,\ldots ,m+n$ and $r=1,2,\ldots.$
We shall need  the following analogues of the families $\Gamma$ and $\Gamma^+$: 
\begin{align*}
&\mathcal{G} =\left\{\beta_{ij}^{(2r-1)}\,:\, (i,j)\in I_0\text{ and }r\geqslant 1\right\} \cup 
\left\{b_{ij}^{( 2r )}\,:\, (i,j)\in I_1\text{ and }r\geqslant 1\right\} \subset\BB,\\
&\mathcal{G}^+=\left\{\beta_{ij}^{(-2r)}\,:\, (i,j)\in I_0\text{ and }r\geqslant 1\right\}\cup
\left\{b_{ij}^{(-2r+1)}\,:\, (i,j)\in I_1\text{ and }r\geqslant 1\right\}\subset\BBd .
\end{align*}
Consider the ascending filtration over $\DBB$ defined by the degree operator $\deg_2$,
$$
\deg_2\ts  \beta_{ij}^{(r)}=r-1 \fand\deg_2\ts  \beta_{ij}^{(-r)}=-r 
$$
for all $i,j=1,\ldots ,m+n$ and $r=1,2,\ldots .$
As before, we write $\bar{\beta}_{ij}^{(r)}$   for the image  of the generator  $ \beta_{ij}^{(r)}$  
in the respective component of the corresponding graded algebra $\gr_2  \DBB$ and we use the notation
$$
\wvr{\mathcal{G}}=\left\{\bar{\beta}_{ij}^{(r)}\,:\,\beta_{ij}^{(r)}\in \mathcal{G}\right\} 
\fand
\wvr{\mathcal{G}}^+=\left\{\bar{\beta}_{ij}^{(-r)}\,:\, \beta_{ij}^{(-r)}\in \mathcal{G}^+\right\}.
$$

By comparing the   relations \eqref{refl1}--\eqref{unidentities} and \eqref{brefl1}--\eqref{bunidentities}, we   conclude  that    the assignments
$\beta_{ij}^{(\pm r)}\mapsto b_{ij}^{(\pm r)}$,
where  $i,j=1,\ldots ,m+n$ and $r\geqslant 1$,
define the algebra epimorphism   
\beq\label{bbmapY}
\DBB\to  \DB.
\eeq
Furthermore, the map given by \eqref{bbmapY} is filtration-preserving and it gives rise to the epimorphism of the corresponding graded algebras
\beq\label{grbbmapY}
\textstyle\gr_2  \DBB\to  \gr_2 \DB .
\eeq
Indeed,   the surjectivity follows by the last assertion of Theorem \ref{reflprro}. 
Our    goal is to show that   \eqref{bbmapY} is an   isomorphism. In order to do so, we shall need the next two lemmas.

\begin{lem}\label{lemabbplus}
 Fix some ordering on  $\wvr{\mathcal{G}}\cup\wvr{\mathcal{G}}^+$. 
Then any element of $\gr_2 \BB$ (resp. $\gr_2 \BBd$) can be written   as a linear combination of ordered monomials in elements of $\wvr{\mathcal{G}}$
(resp. $\wvr{\mathcal{G}}^+$) with at most power 1 for odd generators.
\end{lem} 

\begin{prf}
In parallel with the proof of Theorem \ref{reflprro}, one can employ the first (resp. second)  unitarity identity in \eqref{bunidentities}   to show that the family $\wvr{\mathcal{G}}$ (resp. $\wvr{\mathcal{G}}^+$) generates the algebra $\gr_2 \BB$ (resp. $\gr_2 \BBd$). Hence, any element of $\gr_2 \BB$ (resp. $\gr_2 \BBd$) can be written   as a linear combination of   monomials in the elements of $\wvr{\mathcal{G}}$
(resp. $\wvr{\mathcal{G}}^+$).  From now on we refer to them more briefly as monomials in   $\wvr{\mathcal{G}}$
(resp. $\wvr{\mathcal{G}}^+$). To prove that it is sufficient to consider   ordered monomials only, one needs to employ the corresponding reflection equation \eqref{brefl1} or \eqref{brefl2}. We sketch the proof in the case of $\gr_2 \BBd$. The analogous arguments can be applied to $\gr_2 \BB $ as well.

Clearly, it is sufficient to check that any monomial in  $\wvr{\mathcal{G}}^+$ can be written as a linear combination of ordered monomials in  $\wvr{\mathcal{G}}^+.$ Suppose  $\mu$ is a monomial in $ \mathcal{G}^+$ of the form
$$
\mu =\mu_1\ts  \beta_{ij}^{(-r)}\ts\beta_{kl}^{(-s)}\ts\mu_2
$$
such that $\bar{\beta}_{kl}^{(-s)}$ precedes $\bar{\beta}_{ij}^{(-r)}$ with respect to the chosen ordering (so that $\bar{\mu}$ is not ordered), while $\mu_1$ and $\mu_2$ are some monomials in $ \mathcal{G}^+$. By extracting the coefficients of $u^{r-1}v^{s-1}$ of the matrix entry $e_{ij}\ot e_{kl}$ in \eqref{brefl2} we get
\beq\label{imgaeof}
\beta_{ij}^{(-r)}\ts\beta_{kl}^{(-s)}=
\pm\beta_{kl}^{(-s)}\ts\beta_{ij}^{(-r)}
+\gamma_1
+\gamma_2,
\eeq
where the sign $\pm$ on the right-hand side depends on the parity of the given generators, $\gamma_1$ stands for a linear combination of some elements $\beta_{pq}^{(-t)}$ such that $-t\leqslant -r-s$ and $\gamma_2$ is a linear combination of some monomials of length two, $\beta_{p_1q_1}^{(-t_1)}\beta_{p_2 q_2}^{(-t_2)}$ such that  $-t_1-t_2 < -r-s$. Hence, by taking the image of \eqref{imgaeof} in the $(-r-s)$-component of the corresponding graded algebra, we express $\bar{\mu}$ as a linear combination 
\beq\label{ertz6}
\bar{\mu}= \pm\bar{\mu}_1\ts \bar{\beta}_{kl}^{(-s)}\ts\bar{\beta}_{ij}^{(-r)}\ts \bar{\mu}_2
+  \bar{\mu}_1\ts \bar{\gamma_1}\ts \bar{\mu}_2.
\eeq
Observe that the length of the monomials which appear in the linear combination $\bar{\mu}_1\ts \bar{\gamma_1}\ts \bar{\mu}_2$ is strictly less than the length of $\bar{\mu}$. 
Therefore, we can  continue to apply such a procedure, now starting with the monomials on the right-hand side of \eqref{ertz6}, until, after finitely many steps, we obtain a linear combination of ordered monomials in $\wvr{\mathcal{G}}^+$, as required.
\end{prf}

In the next lemma, we make use of the remaining  reflection equation \eqref{brefl3}. 

\begin{lem}\label{lemauredjaj}
 Fix some ordering on  $\wvr{\mathcal{G}}\cup\wvr{\mathcal{G}}^+$ so that all elements of $\wvr{\mathcal{G}}^+$ precede all elements of $\wvr{\mathcal{G}}$. Then any element of $\gr_2  \DBB$ can be written as a linear combination of ordered monomials in $\wvr{\mathcal{G}}\cup\wvr{\mathcal{G}}^+$ with at most power 1 for odd generators.
\end{lem}

\begin{prf}
It is sufficient to prove that any element  $z\in\gr_2 \DBB$ can be written as a linear combination of   elements of the form $xy$ with $x\in\gr_2\BBd$ and $y\in\gr_2\BB$. Indeed, the assertion of  the lemma then   follows by applying Lemma \ref{lemabbplus} to $x$ and $y$. Moreover, by Lemma  \ref{lemabbplus}, the set $\wvr{\mathcal{G}}\cup\wvr{\mathcal{G}}^+$ generates $\gr_2  \DBB$ so we can assume without loss of generality that $z$ is a monomial in $\wvr{\mathcal{G}}\cup\wvr{\mathcal{G}}^+$. As in the proof of Lemma \ref{lemabbplus}, we shall   describe a procedure  which one can employ to express $z$  in such a way.

 Suppose  $\mu$ is a monomial in $ \mathcal{G} \cup\mathcal{G}^+$ of the form
$$
\mu =\mu_1\ts  \beta_{ij}^{(r)}\ts\beta_{kl}^{(-s)}\ts\mu_2,
$$
where $i,j,k,l=1,\ldots ,m+n$, $r,s=1,2,\ldots$ and
  $\mu_1,\mu_2$ are some monomials in $ \mathcal{G} \cup\mathcal{G}^+$. As $\bar{\beta}_{kl}^{(-s)}$ precedes $\bar{\beta}_{ij}^{(r)}$ with respect to the chosen ordering, the monomial $\bar{\mu}$ is not ordered.
The  $R$-matrices
$
R(u\pm v)
$, which appear in the reflection equation \eqref{brefl3}, are of  the form
	$$
	 R(u\pm v )=1-\frac{P}{u}\sum_{l\geqslant 1}(\mp 1)^l\frac{v^l}{u^{l}}
 .
	$$
For $r,s\geqslant 2$  the degree of the coefficient of $u^{-r}v^{s-1}$ in   $\frac{v^l}{u^l}\beta_{ij}^{kl}(u,v)$
	(resp. $\frac{1}{u}\beta_{ij}^{kl}(u,v)$), where
	$$
	\beta_{ij}^{kl}(u,v)=\left(\beta_{ij}(u)-g_{ij}\right)\left(g_{kl}-\beta_{kl}^+(v)\right),
	$$
	is less than or equal to $r-s-1$ (resp.  $r-s-2$).  
Using these observations and extracting the coefficients of $u^{-r}v^{s-1}$ of the matrix entry $e_{ij}\ot e_{kl}$ in \eqref{brefl3}, we get
\beq\label{imgaeof2}
\beta_{ij}^{(r)}\ts\beta_{kl}^{(-s)}+\gamma_1+\gamma_2
=
\pm\beta_{kl}^{(-s)}\ts\beta_{ij}^{(r)}+\delta_1+\delta_2,
\eeq
where the sign $\pm$ on the right-hand side depends on the parity of the given generators, $\gamma_1,\delta_1$ are  linear combinations of the elements $\beta_{pq}^{(\pm t)}$ of degree less than or equal to $r-s-1$ and $\gamma_2,\delta_2$ are  linear combinations of some  monomials  in $ \mathcal{G} \cup \mathcal{G}^+$ of degree less than or equal to $r-s-2$. Hence, by taking the image of \eqref{imgaeof2} in the $(r-s-1)$-component of the corresponding graded algebra, we express $\bar{\mu}$ as a linear combination
\beq\label{nuredjenimnomii}
\bar{\mu}=\pm \bar{\mu}_1\ts \bar{\beta}_{kl}^{(-s)}\ts\bar{\beta}_{ij}^{(r)}\ts \bar{\mu}_2
+\bar{\mu}_1\ts \bar{\delta}_1\ts \bar{\mu}_2
-\bar{\mu}_1\ts \bar{\gamma_1}\ts \bar{\mu}_2.
\eeq
As with  the proof of Lemma  \ref{lemabbplus}, we can continue to apply such a procedure, now
starting with the monomials on the right-hand side of \eqref{nuredjenimnomii}, until, after finitely many
steps, we obtain a linear combination of   the elements   of the form $xy$ with $x\in\gr_2\BBd$ and $y\in\gr_2\BB$, as required.
\end{prf}

We are now ready to prove the main result of this section.

\begin{thm}
The map \eqref{bbmapY} is an algebra isomorphism.
\end{thm}

\begin{prf}
It is sufficient to check that the  map \eqref{bbmapY} is injective. Therefore, let us assume that its kernel contains some   nonzero element   $x$. Then the image   of   $\bar{x}\in \gr_2 \DBB$ under the map \eqref{grbbmapY}  is trivial, which leads to contradiction. Indeed, by Lemma \ref{lemauredjaj}, the element $\bar{x}$ can be expressed as a linear combination of ordered monomials in $\wvr{\mathcal{G}}\cup\wvr{\mathcal{G}}^+$  with at most power 1 for odd generators. Hence, its image in $\gr_2 \DB $  under the map \eqref{grbbmapY} is a linear combination   of ordered monomials in  $\wvr{\Gamma}\cup \wvr{\Gamma}^+$  with at most power 1 for odd generators. However, by Corollary \ref{koorolar} and the Poincar\'{e}--Birkhoff--Witt Theorem for the enveloping algebra $ \U(\mathfrak{gl}_{m|n}[t,t^{-1}]^{\sigma} )$,  such a linear combination is nonzero.
\end{prf}

\linespread{1.0}

 \section*{Acknowledgement}
This work has been supported in part by Croatian Science Foundation under the project UIP-2019-04-8488.

\end{document}